\documentclass[11pt,a4paper]{art-num}
\usepackage{amsfonts}
\usepackage{amssymb}
\usepackage{amsmath}
\usepackage{hyperref}
\usepackage[page,title]{appendix}
\usepackage{fancyhdr}

\renewcommand{\theequation}{\thesection.\arabic{equation}}
\input{tcilatex}
\reportno{AEI-2009-095}

\begin{document}

\title{Betti numbers of a class of barely $G_{2}$ manifolds}
\author{Sergey Grigorian \\
Max-Planck-Institut f\"{u}r Gravitationsphysik (Albert-Einstein-Institut)\\
Am M\"{u}hlenberg 1\\
D-14476 Golm\\
Germany}
\maketitle

\begin{abstract}
We calculate explicitly the Betti numbers of a class of barely $G_{2}$
manifolds - that is, $G_{2}$ manifolds that are realised as a product of a
Calabi-Yau manifold and a circle, modulo an involution. The particular class
which we consider are those spaces where the Calabi-Yau manifolds are
complete intersections of hypersurfaces in products of complex projective
spaces and the involutions are free acting.
\end{abstract}

\section{Introduction}

One of the key concepts in String and M-theory is the concept of
compactification - here the full $10$- or $11$-dimensional spacetime is
considered to be of the form $M_{4}\times X$ where $M_{4}$ is the
\textquotedblleft large\textquotedblright\ $4$-dimensional visible
spacetime, while $X$ is the \textquotedblleft small\textquotedblright\
compact $6$- or $7$-dimensional Riemannian manifold. Due to considerations
of supersymmetry, these compact manifolds have to satisfy certain conditions
which place restrictions on the geometry. In the case of String theory, the $%
6$-dimensional manifolds have to be Calabi-Yau manifolds - that is K\"{a}%
hler manifolds with vanishing first Chern class. The existence of Ricci-flat
K\"{a}hler metrics for these manifolds has been proven by Yau in 1978 \cite%
{CalabiYau}. One of the first examples of a Calabi-Yau $3$-fold (6 real
dimensions) was the quintic - a degree 5 hypersurface in $\mathbb{CP}^{4}$.
Later, Candelas et al \cite{CandelasCICY1} found the first large class of
Calabi-Yau manifolds - the Complete Intersection Calabi-Yau (CICY)
manifolds, which are given by intersections of hypersurfaces in products of
complex projective spaces. We review the details in section \ref{cicysect}.
Since then even larger classes of Calabi-Yau manifolds have been constructed
- such as Weighted Complete Intersection manifolds \cite{GreeneYAU1}, and
complete intersection manifolds in toric varieties \cite{batyrev-1994}. So
overall there is a very large pool of examples of Calabi-Yau manifolds, and
it is in fact still an open question whether the number of topologically
distinct Calabi-Yau $3$-folds is finite or not. One of the great discoveries
in the study of Calabi-Yau manifolds is Mirror Symmetry \cite%
{Strominger:1996it,MirrorSymBook}. This symmetry first appeared in String
Theory where evidence was found that conformal field theories (CFTs) related
to compactifications on a Calabi-Yau manifold with Hodge numbers $\left(
h_{1,1},h_{2,1}\right) $ are equivalent to CFTs on a Calabi-Yau manifold
with Hodge numbers $\left( h_{2,1},h_{1,1}\right) $. Mirror symmetry is
currently a powerful tool both for calculations in String Theory and in the
study of the Calabi-Yau manifolds and their moduli spaces.

However if we go one dimension higher, and look at compactifications of $M$%
-theory, a natural analogue of a Calabi-Yau manifold in this setting is a $7$%
-dimensional manifold with $G_{2}$ holonomy. These manifolds are also
Ricci-flat, but being odd-dimensional they are real manifolds. The first
examples of $G_{2}$ manifolds have been constructed by Joyce in \cite%
{JoyceG2}. While some work has been done both on the physical aspects of $%
G_{2}$ compactifications (for example \cite%
{Harvey:1999as,Gutowski:2001fm,AcharyaGukov,WittenBeasley} among others) and
on the structure and properties of the moduli space (for example \cite%
{JoyceG2,Lee:2002fa,karigiannis-2007,karigiannis-2007a,GrigorianYau1} among
others), still very little is known about the overall structure of $G_{2}$
moduli spaces. One of the problems is that there are relatively few examples
of $G_{2}$ manifolds, and for the ones that are known it is hard to do any
calculations, because the examples are not very explicit. However there is a
conjectured method of constructing $G_{2}$ manifolds from Calabi-Yau
manifolds, which could potentially yield many new examples of $G_{2}$
manifolds. Here we take a Calabi-Yau $3$-fold $Y$ and let $Z=(Y\times S^{1})/%
\hat{\sigma}$ where $\hat{\sigma}$ acts as antiholomorphic involution on $Y$
and acts as $z\longrightarrow -z$ on the $S^{1}$. In general, the result
will have singularities, and it is still an unresolved question how to
systematically resolve these singularities to obtain a smooth manifold with $%
G_{2}$ holonomy. This construction has been suggested by Joyce in \cite%
{JoyceG2, Joycebook}. A more basic approach is to only consider involutions
without fixed points, so that the resulting manifold $Z$ is smooth.
Manifolds belonging to this class have been called \emph{barely }$G_{2}$ 
\emph{manifolds }in \cite{Harvey:1999as}. Such manifolds do not have the
full $G_{2}$ holonomy, but rather only $SU\left( 3\right) \ltimes \mathbb{Z}%
_{2}$. However, they do share many of the same properties as full $G_{2}$
manifolds, so for many purposes they can play the same role as genuine $%
G_{2} $ manifolds \cite{Harvey:1999as,PartouchePioline}. In particular, if
we consider a specific class of of Calabi-Yau manifolds, such as CICY
manifolds, we can construct a corresponding class of barely $G_{2}$
manifolds rather explicitly. This is what we focus on in this paper. We
first give an overview of $G_{2}$ manifolds and CICY manifolds, and then
describe the algorithm that was used to systematically calculate the Betti
numbers of the barely $G_{2}$ manifolds corresponding to the independent
CICY manifolds.

\textbf{Acknowledgements. }I would like to thank Tristan H\"{u}bsch for the
useful correspondence about CICY Hodge number, and Rahil Baber for the help
with programming.

\section{$G_{2}$ manifolds}

\subsection{\label{g2basicsect}Basics}

We will first review the basics of manifolds with $G_{2}$ holonomy. The $14$%
-dimensional exceptional Lie group $G_{2}\subset SO\left( 7\right) $ is
precisely the group of automorphisms of imaginary octonions, so it preserves
the octonionic structure constants \cite{BaezOcto}. Suppose $x^{1},...,x^{7}$
are coordinates on $\mathbb{R}^{7}$ and let $e^{ijk}=dx^{i}\wedge
dx^{j}\wedge dx^{k}$. Then define $\varphi _{0}$ to be the $3$-form on $%
\mathbb{R}^{7}$ given by 
\begin{equation}
\varphi _{0}=e^{123}+e^{145}+e^{167}+e^{246}-e^{257}-e^{347}-e^{356}.
\label{phi0def}
\end{equation}%
These precisely give the structure constants of the octonions, so $G_{2}$
preserves $\varphi _{0}$. Since $G_{2}$ preserves the standard Euclidean
metric $g_{0}$ on $\mathbb{R}^{7}$, it preserves the Hodge star, and hence
the dual $4$-form $\ast \varphi _{0},$ which is given by 
\begin{equation}
\ast \varphi
_{0}=e^{4567}+e^{2367}+e^{2345}+e^{1357}-e^{1346}-e^{1256}-e^{1247}.
\label{sphi0def}
\end{equation}

Now suppose $X$ is a smooth, oriented $7$-dimensional manifold. A $G_{2}$
structure $Q$ on $X$ is a principal subbundle of the frame bundle $F$, with
fibre $G_{2}$. However we can also uniquely define $Q$ via $3$-forms on $X.$
Define a $3$-form $\varphi $ to be \emph{positive }if we locally can choose
coordinates such that $\varphi $ is written in the form (\ref{phi0def}) -
that is for every $p\in X$ there is an isomorphism between $T_{p}X$ and $%
\mathbb{R}^{7}$ such that $\left. \varphi \right\vert _{p}=\varphi _{0}$.
Using this isomorphism, to each positive $\varphi $ we can associate a
metric $g$ and a Hodge dual $\ast \varphi $ which are identified with $g_{0}$
and $\ast \varphi _{0}$ under this isomorphism. It is shown in \cite%
{Joycebook} that there is a $1-1$ correspondence between positive $3$-forms $%
\varphi $ and $G_{2}$ structures $Q$ on $X$.

So given a positive $3$-form $\varphi $ on $X$, it is possible to define a
metric $g$ associated to $\varphi $ and this metric then defines the Hodge
star, which in turn gives the $4$-form $\ast \varphi $. Thus although $\ast
\varphi $ looks linear in $\varphi $, it actually is not, so sometimes we
will write $\psi =\ast \varphi $ to emphasize that the relation between $%
\varphi $ and $\ast \varphi $ is very non-trivial.

It turns out that the holonomy group $Hol\left( X,g\right) \subseteq G_{2}$
if and only if $X$ has a torsion-free $G_{2}$ structure \cite{Joycebook}. In
this case, the invariant $3$-form $\varphi $ satisfies%
\begin{equation}
d\varphi =d\ast \varphi =0  \label{torsionfreedef}
\end{equation}%
and equivalently, $\nabla \varphi =0$ where $\nabla $ is the Levi-Civita
connection of $g$. So in fact, in this case $\varphi $ is harmonic.
Moreover, if $Hol\left( X,g\right) \subseteq G_{2}$, then $X$ is Ricci-flat.
The holonomy group is precisely $G_{2}$ only if the first Betti number $%
b_{1} $ vanishes.

Special holonomy manifolds play a very important role in string and $M$%
-theory because of their relation to supersymmetry. In general, if we
compactify string or $M$-theory on a manifold of special holonomy $X$ the
preservation of supersymmetry is related to existence of covariantly
constant spinors (also known as parallel spinors). In fact, if all bosonic
fields except the metric are set to zero, and a supersymmetric vacuum
solution is sought, then in both string and $M$-theory, this gives precisely
the equation 
\begin{equation}
\nabla \xi =0  \label{covconstspinor}
\end{equation}%
for a spinor $\xi $. As lucidly explained in \cite{AcharyaGukov}, condition (%
\ref{covconstspinor}) on a spinor immediately implies special holonomy. Here 
$\xi $ is invariant under parallel transport, and is hence invariant under
the action of the holonomy group $Hol\left( X,g\right) $. This shows that
the spinor representation of $Hol\left( X,g\right) $ must contain the
trivial representation. For $Hol\left( X,g\right) =SO\left( n\right) $, this
is not possible since the spinor representation is reducible, so $Hol\left(
X,g\right) \subset SO\left( n\right) $. In particular, Calabi-Yau 3-folds
with $SU\left( 3\right) $ holonomy admit two covariantly constant spinors
and $G_{2}$ holonomy manifolds admit only one covariantly constant spinor.
Hence eleven-dimensional supergravity compactified on a $G_{2}$ holonomy
manifold gives rise to a $\mathcal{N}=1$ effective theory. From \cite%
{AcharyaGukov},\cite{WittenBeasley} and \cite{Gutowski:2001fm} we know that
the deformations of the $G_{2}$ $3$-form $\varphi $ give $b_{3}$ real moduli
which combine with the deformations of the supergravity $3$-form $C$ to give 
$b_{3}$ complex moduli. Together with modes of the gravitino, this gives $%
b_{3}$ chiral multiplets. Decomposition of the $C$-field also gives $b_{2}$
abelian gauge fields, which again combine with gravitino modes to give $%
b_{2} $ vector multiplets. The structure of the moduli space has been
studied in detail in \cite{GrigorianYau1}.

Examples of compact $G_{2}$ manifolds have been first constructed by Joyce 
\cite{JoyceG2} as orbifolds $T^{7}/\Gamma $ for a discrete group $\Gamma $.
There $\Gamma $ is taken to be a finite group of diffeomorphisms of $T^{7}$
preserving the flat $G_{2}$-structure on $T^{7}$. The resulting orbifold
will have a singular set coming from the fixed point of the action of $%
\Gamma $, and these singularities are resolved by gluing ALE\ spaces with
holonomy $SU\left( 2\right) $ or $SU\left( 3\right) $.

\subsection{$G_{2}$ manifolds from Calabi-Yau manifolds}

A simple way to construct a manifold with a torsion-free $G_{2}$ structure
is to consider $X=Y\times S^{1}$ where $Y$ is a Calabi-Yau $3$-fold. Define
the metric and a $3$-form on $X$ as 
\begin{eqnarray}
g_{X} &=&d\theta ^{2}\times g_{Y}  \label{metCY} \\
\varphi &=&d\theta \wedge \omega +\func{Re}\Omega  \label{phiCY}
\end{eqnarray}%
where $\theta $ is the coordinate on $S^{1},$ $\omega $ is the K\"{a}hler
form on $Y$ and $\Omega $ is the holomorphic $3$-form on $Y$. This then
defines a torsion-free $G_{2}$ structure, with 
\begin{equation}
\ast \varphi =\frac{1}{2}\omega \wedge \omega -d\theta \wedge \func{Im}%
\Omega .  \label{psiCY}
\end{equation}%
However, the holonomy of $X$ in this case is $SU\left( 3\right) \subset
G_{2} $. From the K\"{u}nneth formula we get the following relations between
the Betti numbers of $X$ and the Hodge numbers of $Y$ 
\begin{eqnarray*}
b_{1} &=&1\ \ \  \\
b_{2} &=&h_{1,1} \\
b_{3} &=&h_{1,1}+2\left( h_{2,1}+1\right) \ \text{\ }
\end{eqnarray*}%
In \cite{JoyceG2} and \cite{Joycebook}, Joyce describes a possible
construction of a smooth manifold with holonomy equal to $G_{2}$ from a
Calabi-Yau manifold $Y$. So suppose $Y$ is a Calabi-Yau $3$-fold as above.
Then suppose $\sigma :Y\longrightarrow Y$ is an antiholomorphic isometric
involution on $Y$, that is, $\chi $ preserves the metric on $Y$ and
satisfies 
\begin{subequations}
\label{involprop}
\begin{eqnarray}
\sigma ^{2} &=&1 \\
\sigma ^{\ast }\left( \omega \right) &=&-\omega \\
\sigma ^{\ast }\left( \Omega \right) &=&\bar{\Omega}.
\end{eqnarray}%
Such an involution $\sigma $ is known as a \emph{real structure }on $Y$.
Define now a quotient given by 
\end{subequations}
\begin{equation}
Z=\left( Y\times S^{1}\right) /\hat{\sigma}  \label{barelydefine}
\end{equation}%
where $\hat{\sigma}$:$Y\times S^{1}\longrightarrow Y\times S^{1}$ is defined
by $\hat{\sigma}\left( y,\theta \right) =\left( \sigma \left( y\right)
,-\theta \right) $. The $3$-form $\varphi $ defined on $Y\times S^{1}$ by (%
\ref{phiCY}) is invariant under the action of $\hat{\sigma}$ and hence
provides $Z$ with a $G_{2}$ structure. Similarly, the dual $4$-form $\ast
\varphi $ given by (\ref{psiCY}) is also invariant. Generically, the action
of $\sigma $ on $Y$ will have a non-empty fixed point set $N$, which is in
fact a special Lagrangian submanifold on $Y$ \cite{Joycebook}. This gives
rise to orbifold singularities on $Z$. The singular set is two copies of $Z$%
. It is conjectured that it is possible to resolve each singular point using
an ALE $4$-manifold with holonomy $SU\left( 2\right) $ in order to obtain a
smooth manifold with holonomy $G_{2}$, however the precise details of the
resolution of these singularities are not known yet. We will therefore
consider only free-acting involutions, that is those without fixed points.

Manifolds defined by (\ref{barelydefine}) with a freely acting involution
were called \emph{barely }$G_{2}$ \emph{manifolds }by Harvey and Moore in 
\cite{Harvey:1999as}. The cohomology of barely $G_{2}$ manifolds is
expressed in terms of the cohomology of the underlying Calabi-Yau manifold $%
Y $: 
\begin{subequations}
\label{barelycoho}
\begin{eqnarray}
H^{2}\left( Z\right) &=&H^{2}\left( Y\right) ^{+} \\
H^{3}\left( Z\right) &=&H^{2}\left( Y\right) ^{-}\oplus H^{3}\left( Y\right)
^{+}
\end{eqnarray}%
Here the superscripts $\pm $ refer to the $\pm $ eigenspaces of $\sigma
^{\ast }$. Thus $H^{2}\left( Y\right) ^{+}$ refers to two-forms on $Y$ which
are invariant under the action of involution $\sigma $ and correspondingly $%
H^{2}\left( Y\right) ^{-}$ refers to two-forms which are odd under $\sigma $%
. Wedging an odd two-form on $Y$ with $d\theta $ gives an invariant $3$-form
on $Y\times S^{1}$, and hence these forms, together with the invariant $3$%
-forms $H^{3}\left( Y\right) ^{+}$ on $Y$, give the three-forms on the
quotient space $Z$. Also note that $H^{1}\left( Z\right) $ vanishes, since
the $1$-form on $S^{1}$ is odd under $\hat{\sigma}$. Now, given a $3$-form
on $Y$, its real part will be invariant under $\sigma $, hence $H^{3}\left(
Y\right) ^{+}$ is essentially the real part of $H^{3}\left( Y\right) $.
Therefore the Betti numbers of $Z$ in terms of Hodge numbers of $Y$ are 
\end{subequations}
\begin{subequations}
\label{barelybetti}
\begin{eqnarray}
b^{1} &=&0 \\
b^{2} &=&h_{1,1}^{+} \\
b^{3} &=&h_{1,1}^{-}+h_{2,1}+1
\end{eqnarray}%
Hence in order to construct barely $G_{2}$ manifolds we need to be able to
find involutions of Calabi-Yau manifolds and determine the action of the
involution on $H^{1,1}\left( Y\right) $. A relatively large class of
Calabi-Yau manifolds for which this is not hard to do are the complete
intersection Calabi-Yau manifolds. We review the properties of these
manifolds in the next section.

\section{Complete intersection Calabi-Yau manifolds}

\subsection{Basics}

\label{cicysect}Complete intersection Calabi-Yau (CICY) manifolds were the
first major class of Calabi-Yau manifolds which was discovered by Candelas
et al in \cite{CandelasCICY1}. Such a manifold $M$ is defined as a complete
intersection of $K$ hypersurfaces in a product of $m$ complex projective
spaces $W=\mathbb{CP}^{n_{1}}\times ...\times \mathbb{CP}^{n_{m}}$. Each
hypersurface is defined as the zero set of a homogeneous holomorphic
polynomial 
\end{subequations}
\begin{equation}
f^{a}\left( z_{\ r}^{\mu }\right) =0\ \ a=1,...,K.  \label{homopoly}
\end{equation}%
Each such polynomial is homogeneous of degree $q_{a}^{r}$ with respect to
the homogeneous coordinates of $\mathbb{CP}^{n_{r}}$. By complete
intersection it is meant that the $K$-form 
\begin{equation*}
\Theta =df^{1}\wedge ...\wedge df^{K}
\end{equation*}%
does not vanish on $M$. This condition ensures that the resulting manifold
is defined globally. In order for $M$ to be a $3$-fold, we obviously need 
\begin{equation}
K=\sum_{i=1}^{m}n_{i}-3.  \label{kcond}
\end{equation}%
The standard notation for a CICY manifold is a $m\times (K+1)$ array of the
form 
\begin{equation}
\left[ \left. n\right\Vert q\right]   \label{cicyarray}
\end{equation}%
where $n$ is a column $m$-vector whose entries $n_{r}$ are the dimensions of
the $\mathbb{CP}^{n_{r}}$ factors, and $q$ is a $m\times K$ matrix with
entries $q_{a}^{r}$ which give the degrees of the polynomials in the
coordinates of each of the $\mathbb{CP}^{n_{r}}$ factor. Each such array
defining a CICY is known as a \emph{configuration} \emph{matrix}, while an
equivalence class of configuration matrices under permutation of all rows
and all columns belonging to $q$ is called a \emph{configuration}. Clearly
each such a permutation defines exactly the same manifold.

As it was shown in \cite{CandelasCICY1}, Chern classes can be computed
directly from the defining quantities $n$ and $q$. In particular, we
immediately get the condition for a vanishing first Chern class:%
\begin{equation}
n_{r}+1=\sum_{a=1}^{K}q_{a}^{r}\ \ \ \ \forall r  \label{1stcherncond}
\end{equation}%
That is, the sum of entries of in each row of $q$ must equal to the
dimension of the corresponding $\mathbb{CP}^{n_{r}}$ factors. This is hence
precisely the condition for the complete intersection manifold to be
Calabi-Yau. Moreover from the expressions for Chern classes, an expression
for the Euler number is also obtained. This is given by 
\begin{equation}
\chi _{E}\left( M\right) =\left[ \left(
\sum_{r,s,t=1}^{m}c_{3}^{rst}x_{r}x_{s}x_{t}\right) \cdot
\prod_{b=1}^{K}\left( \sum_{u=1}^{m}q_{b}^{u}x_{u}\right) \right] _{\text{%
coefficient of }\prod_{r=1}^{m}\left( x_{r}\right) ^{n_{r}}}
\label{EulerChar}
\end{equation}%
where 
\begin{equation*}
c_{3}^{rst}=\frac{1}{3}\left( \left( n_{r}+1\right) \delta
^{rst}-\sum_{a=1}^{K}q_{a}^{r}q_{a}^{s}q_{a}^{t}\right)
\end{equation*}%
and $\delta ^{rst}$\thinspace $=1$ for $r=s=t$ and vanishes otherwise.

Varying the coefficients of polynomials in a CICY\ configuration generally
corresponds to complex structure deformations, but as it was shown in \cite%
{Green:1987rw}, there is no one to one correspondence. So it is said that
each configuration corresponds to a partial deformation class. There are
also various identities which relate different configurations, so not all
configurations are independent. There are however 7868 independent
configurations. A method for calculating Hodge numbers of the CICY manifolds
has been found by Green and H\"{u}bsch in \cite{Green:1987rw} and in \cite%
{GreenAllHodge:1987cr} Green, H\"{u}bsch and L\"{u}tken calculated the Hodge
numbers for each of the 7868 configurations. They found there to be 265
unique pairs of Hodge numbers. Unfortunately, the original data with the
CICY Hodge numbers has been lost, and the original computer code by H\"{u}%
bsch has been written in a curious mix of $C$ and $Pascal$ so the original
code had to be rewritten in standard $C$ in order to be able to recompile
the list of Hodge numbers for CICY manifolds, which is necessary to be able
to calculate the Betti numbers of corresponding barely $G_{2}$ manifolds.

\subsection{Involutions}

Antiholomorphic involutions of projective spaces have been classified in 
\cite{PartouchePioline}, and here we briefly review their results. First
consider involutions of a single projective space $\mathbb{CP}^{n}$. Suppose
we have homogeneous coordinates $\left( z_{0},z_{1},...,z_{n}\right) $ on $%
\mathbb{CP}^{n}$, then we can represent an anti-holomorphic involution $%
\sigma $ by a matrix $M$ which acts as 
\begin{equation}
z_{i}\longrightarrow M_{ij}\bar{z}_{j}  \label{involmatrix}
\end{equation}%
Without loss of generality we fix $\det M=1$ since multiplication by any
non-zero complex number still gives the same involution. Moreover,
involutions which differ only by a holomorphic change of basis can be
regarded to be the same.

Also $\sigma ^{2}=1$ must be true projectively, so we get 
\begin{equation}
M\bar{M}=\lambda I\text{.}  \label{mmbarinvol}
\end{equation}%
Taking the determinant of (\ref{mmbarinvol}) we find that $\lambda ^{n+1}=1$%
, and taking the trace we see that $\lambda $ is real. Thus $\lambda =1$ for 
$n$ even and $\lambda =\pm 1$ for $n$ odd. The involution $\sigma $ is
required to be an isometry - that is, it must preserve the standard
Fubini-Study metric of $\mathbb{CP}^{n}.$ Together with previous
restrictions on $M,$ this gives the condition%
\begin{equation}
MM^{\dag }=I.  \label{mmdaginvol}
\end{equation}%
Combining (\ref{mmbarinvol}) and (\ref{mmdaginvol}), we see that for $%
\lambda =1$ these equations imply that $M$ is symmetric, and for $\lambda
=-1 $ that $M$ is antisymmetric. Moreover, due to (\ref{mmbarinvol}), the
real and imaginary parts of $M$ commute, and so can be simultaneously
brought into a canonical form - diagonal for $\lambda =1$ and block-diagonal
for $\lambda =-1$. Another change of basis can be used to normalize the
coefficients. Hence we get two distinct antiholomorphic involutions 
\begin{subequations}
\begin{eqnarray}
A &:&\left( z_{0},z_{1},...,z_{n}\right) \longrightarrow \left( \bar{z}_{0},%
\bar{z}_{1},...,\bar{z}_{n}\right) \\
B &:&\left( z_{0},z_{1},...,z_{n-1},z_{n}\right) \longrightarrow \left( -%
\bar{z}_{1},\bar{z}_{0},...,-\bar{z}_{n},\bar{z}_{n-1}\right) .
\end{eqnarray}%
The involution $A$ corresponds to $\lambda =+1$ and is defined for $n$ both
odd and even, whereas the involution $B$ corresponds to $\lambda =-1$ and is
only defined for $n$ odd. An important difference between the two
involutions is that $A$ has a fixed point set $\left\{ z_{i}=\bar{z}%
_{i}\right\} $, whereas $B$ acts freely without any fixed points.

So far we considered antiholomorphic involutions of a single projective
space, but in general we are interested in products of projective spaces, so
we should also consider involutions which mix different factors. As pointed
out in \cite{PartouchePioline}, the only possibility for this is two
exchange two identical projective factors $\mathbb{CP}^{n}$, giving another
involution $C$: 
\end{subequations}
\begin{equation}
C:\left( \left\{ y_{i}\right\} ;\left\{ z_{i}\right\} \right)
\longrightarrow \left( \left\{ \bar{z}_{i}\right\} ;\left\{ \bar{y}%
_{i}\right\} \right) .  \label{Cinvol}
\end{equation}%
This involution clearly has a fixed point set $\left\{ y_{i}=\bar{z}%
_{i}\right\} $.

Now that we have antiholomorphic involutions of projective spaces, we can
use these to construct barely $G_{2}$ manifolds from CICY manifolds, as in (%
\ref{barelydefine}). In general we must either have an involution acting on
each projective factor - either involutions $A$ or $B$ on single factors or
involution $C$ on a pair of identical projective factors.

Given a CICY configuration matrix, we will denote the resulting barely $%
G_{2} $ manifold by the same configuration matrix, but indicating in the
first column of the configuration matrix which involutions are acting on
each projective factor. These actions will be denoted by $\bar{n}$, $\hat{n}$
and $%
\begin{array}{c}
\overset{\frown }{n} \\ 
\underset{\smile }{n}%
\end{array}%
$ for involutions $A$, $B$ and $C$, respectively. For example, consider the
configuration matrix:

\begin{equation}
\left[ 
\begin{array}{c}
\widehat{1} \\ 
\overset{\frown }{1} \\ 
\underset{\smile }{1} \\ 
\overline{2} \\ 
\overline{3}%
\end{array}%
\right\Vert \left. 
\begin{array}{ccccc}
0 & 0 & 0 & 0 & 2 \\ 
0 & 0 & 1 & 1 & 0 \\ 
0 & 0 & 1 & 1 & 0 \\ 
1 & 1 & 1 & 0 & 0 \\ 
1 & 1 & 0 & 1 & 1%
\end{array}%
\right] ^{1,39}  \label{configex1}
\end{equation}%
This denotes the barely $G_{2}$ manifolds constructed from CICY with the
same configuration matrix but with involution $A$ acting on the $\mathbb{CP}%
^{2}$ and $\mathbb{CP}^{3}$ factors, involution $B$ acting on the first
remaining $\mathbb{CP}^{1}$ factor and involution $C$ acting on the
remaining $\mathbb{CP}^{1}$ $\times \mathbb{CP}^{1}$. The superscripts $%
\left( 1,39\right) $ give the Betti numbers $b^{2}$ and $b^{3}$ of the
resulting $7$-manifold. Note that since this example includes the action of
involution $B$ which has no fixed points, the full involution acting on the
whole CICY is also free, so the resulting space is a smooth barely $G_{2}$
manifold.

When the projective space involution restricts to the complete intersection
space, conditions are imposed on the coefficients of the defining
homogeneous equations. Thus the involutions must be compatible with the
defining equations, and this may not always be possible. In particular, the
invariance of the defining equations under the involution implies that the
transformed equations must be equivalent to the original equations. Let us
use the configuration matrix (\ref{configex1}) to demonstrate this. Let $%
u_{i}$, $v_{i}$, $w_{i}$ for $i=0,1$ be the homogeneous coordinates on the $%
\mathbb{CP}^{1}$ spaces, let $y_{j}$ for $j=0,1,2$ be coordinates on $%
\mathbb{CP}^{2}$ and $z_{k}$ for $k=0,1,2,3$ be the homogeneous coordinates
on the $\mathbb{CP}^{3}$ factor. Then the original defining equations are 
\begin{equation}
\left\{ 
\begin{array}{c}
f_{1}\left( y,z\right) =f_{2}\left( y,z\right) =0 \\ 
g_{1}\left( v,w,y\right) =g_{2}\left( v,w,z\right) =0 \\ 
h\left( u,z\right) =0%
\end{array}%
\right. ,  \label{ex1defeq}
\end{equation}%
where the $f_{i}$ and $g_{i}$ are polynomials homogeneous of degree $1$ in
their variable and $h$ is a polynomial which is homogeneous of degree $2$ in 
$u_{i}$ and of degree $1$ in $z_{k}$. Under the involution presented in (\ref%
{configex1}), after taking the complex conjugates, these equations become 
\begin{equation}
\left\{ 
\begin{array}{c}
\bar{f}_{1}\left( y,z\right) =\bar{f}_{2}\left( y,z\right) =0 \\ 
\bar{g}_{1}\left( w,v,y\right) =\bar{g}_{2}\left( w,v,z\right) =0 \\ 
\bar{h}\left( \hat{u},z\right) =0%
\end{array}%
\right. ,
\end{equation}%
where $\hat{u}_{2k}=-u_{2k+1}$ and $\hat{u}_{2k+1}=u_{2k}$. Then for some
complex numbers $\lambda _{1},\lambda _{2}$ and $\lambda _{3}$ we must have 
\begin{subequations}
\begin{eqnarray}
g_{1}\left( v,w,y\right) &=&\lambda _{1}\bar{g}_{1}\left( w,v,y\right) \ \ \ 
\label{consist1} \\
g_{2}\left( v,w,z\right) &=&\lambda _{2}\bar{g}_{2}\left( w,v,z\right)
\label{consist2} \\
h\left( u,z\right) &=&\lambda _{3}\bar{h}\left( \hat{u},z\right)
\label{consist3}
\end{eqnarray}%
and for some matrix $M$ in $GL\left( 2,\mathbb{C}\right) $ we must have 
\end{subequations}
\begin{equation}
\ \text{and \ \ \ }\left( 
\begin{array}{c}
f_{1}\left( y,z\right) \\ 
f_{2}\left( y,z\right)%
\end{array}%
\right) =M\left( 
\begin{array}{c}
\bar{f}_{1}\left( y,z\right) \\ 
\bar{f}_{2}\left( y,z\right)%
\end{array}%
\right) .  \label{consist4}
\end{equation}%
For consistency in (\ref{consist1}) and (\ref{consist2}), we find that $%
\lambda _{1}\bar{\lambda}_{1}=1$ and $\lambda _{2}\bar{\lambda}_{2}=1$.
Without loss of generality, we can set $\lambda _{1}=\lambda _{2}=1$. From (%
\ref{consist3}), we have 
\begin{equation}
h\left( u,z\right) =\lambda _{3}\bar{h}\left( \hat{u},z\right) =\lambda _{3}%
\bar{\lambda}_{3}h\left( {\Hat {\Hat u}},z\right) =\lambda _{3}\bar{\lambda}%
_{3}h\left( u,z\right) .  \label{consist3a}
\end{equation}%
Here we have used the fact that $h\left( u,z\right) $ is of degree $2$ in $%
u_{i}$, so even though $\hat{\hat{u}}=-u$, the minus sign cancels, and we
get $\lambda _{3}\bar{\lambda}_{3}=1$. So we can set $\lambda _{3}=1$
without loss of generality. In order for (\ref{consist4}) to be consistent,
we find that we must have $M\bar{M}=I,$ but $M=I$ satisfies this condition
and so fulfills the consistency criteria. We can see that all these
conditions on the coefficients of the defining polynomials halve the number
of possible choices for the coefficients. This also shows that not all
combinations of involutions are possible. In particular, suppose if we
wanted a $B$ involution to act on the $\mathbb{CP}^{3}$ factor. Then since $%
\hat{\hat{z}}=-z,$ and $h\left( u,z\right) $ is of degree $1$ in $z$, from (%
\ref{consist3a}) we would get that $\lambda _{3}\bar{\lambda}_{3}=-1,$ which
is clearly not possible. Also, the $C$ involution is not always possible -
the configuration must be invariant under the interchange of factors.

In order to construct all possible barely $G_{2}$ manifolds from CICY
manifolds, we must be able to find all possible involutions of a given CICY
configuration. Since we want freely acting involutions, we only consider
those combinations of involutions which contain a $B$ involution.

The overall strategy is the following. We first find all possible
combinations of $C$ involutions, and then for each such combination we find
the possible $B$ involutions. The remaining factors which do not have any
involutions acting on them get an $A$ involution.

Suppose we have a configuration matrix with $m$ rows and $K$ columns - that
is we have $K$ hypersurfaces in a product of $m$ projective factors. Let the
coordinates be labelled by $x^{1},...,x^{m}$ and let the homogeneous
polynomials be $f_{1},...,f_{K}$. So the intersection of hypersurfaces is
given by 
\begin{equation}
f_{1}=f_{2}=...=f_{K}=0  \label{hypefcond}
\end{equation}

We want to check whether a $C$ involution is possible on the first two
factors. For this we assume that the two factors are of the same dimension,
as this is a basic necessary condition for a $C$ involution. Then we have to
make sure that after the interchange of $x^{1}$ and $x^{2}$ the new set of
homogeneous equations is equivalent to (\ref{hypefcond}). This is true if
and only if under the interchange of $x^{1}$ and $x^{2}$ the polynomials
remain the same up to a change of ordering. In terms of the configuration
matrix this means that under the interchange of two rows the matrix remains
invariant up to a permutation of the columns. For more than one $C$
involution acting on the same configuration matrix, we thus require that
under the full set of row interchanges the matrix remains invariant up to a
permutation of the columns.

To find all the possible $C$ involutions for a given configuration matrix we
do an exhaustive search of all possibilities. First we find all the possible
combinations of pairs of rows that correspond to projective factors of equal
dimensions. Then for each such combination of pairs we check if under the
interchange of rows in each pair the configuration matrix stays invariant up
to a reordering of columns. If this is true, then it is possible to have $C$
involutions acting on each of these pairs of rows. This procedure then gives
us the full set $\mathcal{C}=\left\{ C_{1},...,C_{N}\right\} $ of all
possible combinations of $C$ involutions acting on the configuration matrix.

Now given all the possible $C$ involutions on a configuration matrix, for
each such combination $C_{i}\in \mathcal{C}$ we need to find the possible $B$
involutions. Suppose we have a configuration matrix as before, and we want
to check whether a $B$ involution is possible on the first projective
factor. The basic necessary condition is that the dimension of this
projective factor is odd. Then we need to make sure that the new set of
homogeneous equations is equivalent to the old set. Let $\mathcal{I}$ be the
set of columns which have non-zero entries in the first row - or
equivalently, the set of polynomials that involve $x^{1}$. First suppose
that all columns in $\mathcal{I}$ are distinct. Then for each $i\in \mathcal{%
I}$ we require 
\begin{equation}
f_{i}\left( z^{1},...\right) =\lambda _{i}\bar{f}_{i}\left( \hat{z}%
^{1},...\right)
\end{equation}%
for some constant $\lambda _{i}\in \mathbb{C}$. As in (\ref{consist3a}), we
then have the consistency requirement%
\begin{equation}
f_{i}\left( z^{1},...\right) =\lambda _{i}\bar{f}_{i}\left( \hat{z}%
^{1},...\right) =\lambda _{i}\bar{\lambda}_{i}f_{i}\left( \hat{\hat{z}}%
^{1},...\right)  \label{hypefcond2}
\end{equation}%
However, $\hat{\hat{z}}^{1}=-z^{1},$ but $f_{i}$ is homogeneous of degree $%
q_{\ i}^{1}$ in $z^{1}$, so $f_{i}\left( \hat{\hat{z}}^{1},...\right)
=\left( -1\right) ^{q_{i}^{1}}f_{i}\left( z^{1},...\right) .$ Hence in order
for (\ref{hypefcond2}) to be consistent, $q_{i}^{1}$ needs to be even for
each $i$. If this is true, then we can have a $B$ involution on the first
projective factor.

More generally, however, suppose that we have some identical columns in $%
\mathcal{I}.$ In particular assume that columns $k_{1},...,k_{r}\in \mathcal{%
I}$ are all identical, and that the remaining columns in $\mathcal{I}$ are
distinct from these. These columns correspond to polynomials which have the
same degrees in projective space coordinates. We can have an involution $B$
if and only if 
\begin{equation*}
f_{k_{1}}=f_{k_{2}}=...=f_{k_{r}}=0\Longleftrightarrow \hat{f}_{k_{1}}=\hat{f%
}_{k_{2}}=...=\hat{f}_{k_{r}}=0.
\end{equation*}%
So for some matrix $M\in GL\left( r,\mathbb{C}\right) $ we must have \ \ 
\begin{equation}
\left( 
\begin{array}{c}
f_{k_{1}}\left( z^{1},...\right) \\ 
... \\ 
f_{k_{r}}\left( z^{1},...\right)%
\end{array}%
\right) =M\left( 
\begin{array}{c}
\bar{f}_{k_{1}}\left( \hat{z}^{1},...\right) \\ 
... \\ 
\bar{f}_{k_{r}}\left( \hat{z}^{1},...\right)%
\end{array}%
\right) .  \label{hypefcond3}
\end{equation}%
From (\ref{hypefcond3}) we have the consistency condition 
\begin{equation}
\left( 
\begin{array}{c}
f_{k_{1}}\left( z^{1},...\right) \\ 
... \\ 
f_{k_{r}}\left( z^{1},...\right)%
\end{array}%
\right) =M\bar{M}\left( 
\begin{array}{c}
f_{k_{1}}\left( \hat{\hat{z}}^{1},...\right) \\ 
... \\ 
f_{k_{r}}\left( \hat{\hat{z}}^{1},...\right)%
\end{array}%
\right) =\left( -1\right) ^{Q}M\bar{M}\left( 
\begin{array}{c}
f_{k_{1}}\left( z^{1},...\right) \\ 
... \\ 
f_{k_{r}}\left( z^{1},...\right)%
\end{array}%
\right) ,  \label{hypefcond3a}
\end{equation}%
where $Q=q_{k_{1}}^{1}+...+q_{k_{r}}^{1}.$ If $r$ is even, then we can
always find a block-diagonal real matrix $M$ such that $M\bar{M}=M^{2}=-I$,
so in this case the condition (\ref{hypefcond3a}) is always consistent,
independent of the parity of $Q$. For example for $r=2$ we could set $%
M=\left( 
\begin{array}{cc}
0 & 1 \\ 
-1 & 0%
\end{array}%
\right) $. However if $r$ is odd, then it is not possible to find a matrix
which satisfies $M\bar{M}=-I$, so we then cannot have $Q$ odd.

To find all possible $B$ involutions, we again proceed with an exhaustive
search. We look for all possible combinations of $B$ involutions for each
combination of $C$ involutions $C_{i}\in \mathcal{C}.$ First we find the set 
$\mathcal{R}$ of all possible combinations of rows such that the dimensions
of the corresponding projective factors are odd, and such that these rows do
not have a $C$ involution from $C_{i}$ acting on them. Given a combination $%
R\in \mathcal{R},$ we want to check if it is possible to have a $B$
involution acting on each row in $R$. We look for the set $\mathcal{I}$ of
columns which have a non-zero entry in at least one of the rows in $R$. The
set $\mathcal{I}$ is then split into maximal subsets of identical columns.
For each such subset we evaluate $Q$ as above, and if for some subset of
size $r$ $rQ$ is odd, then the consistency condition (\ref{hypefcond3a}) is
not fulfilled, and so the combination of rows $R$ does not admit a $B$
involution.

The above algorithm has been implement in the programming language C. After
running the algorithm, for each configuration matrix in the original list of
7868 CICY configurations we find the possible combinations of $C$%
-involutions, and for each combination of $C$-involution all the possible
combinations of $B$ involutions. Since we are interested in manifolds with
free-acting involutions, we are only concerned with those configuration that
admit a $B$-involution. It turns out that a total of 4652 configurations do
admit a $B$-involution, out of which 153 have unique pairs of Hodge numbers.
The Hodge pairs for which there exist configurations that admit a $B$
involutions are listed in (\ref{Binvolpairs})%
\begin{equation}
\begin{tabular}{ll}
$h_{1,1}$ & $h_{2,1}$ \\ 
$1$ & $65,73,89$ \\ 
$2$ & $50+2k\ $for $k=0,...,13,18$ \\ 
$3$ & $31+2k\ $for $k=0,2,3,...,17,19,22$ \\ 
$4$ & $26+2k$ \ for $k=0,1,...,19,21$ \\ 
$5$ & \thinspace $25+2k$ $\ $for $k=0,1,...,18$ \\ 
$6$ & $24+2k$ \ for $k=0,1,...,13,15$ \\ 
$7$ & $23+2k$ for $k=0,1,...,10,12,13$ \\ 
$8$ & $22+2k$ for $k=0,...,11$ \\ 
$9$ & $21+2k$ for $k=0,...,9$ \\ 
$10$ & $20+2k$ for $k=0,...,7$ \\ 
$11$ & $19+2k$ for $k=0,...,6$ \\ 
$12$ & $18+2k$ for $k=0,...,3,5$ \\ 
$13$ & $17+2k$ for $k=0,...,4$ \\ 
$14$ & $16+2k$ for $k=0,1,3$ \\ 
$15$ & $15,21$ \\ 
$16$ & $20$ \\ 
$19$ & $19$%
\end{tabular}
\label{Binvolpairs}
\end{equation}%
As we can see there is a clear pattern - all these pairs of Hodge numbers
have an even sum. In fact the only pairs of Hodge number that have an even
sum but do not admit any $B$ involutions are $\left( 2,46\right) ,$ $\left(
2,64\right) ,\left( 3,27\right) \,$and $\left( 3,33\right) $ .

\section{Barely $G_{2}$ manifolds}

\subsection{Betti numbers}

Now that we have found the CICY involutions, we can calculate the Betti
numbers of the corresponding barely $G_{2}$ manifolds. Thus we need to find
the harmonic forms on these manifolds. As we know from section \ref%
{g2basicsect}, for this we only to determine the stabilizer of the
involution $\sigma $ acting on the $H^{1,1}\left( Y\right) $ of a CICY
manifold $Y$. Suppose $h_{1,1}=m,$ the number of complex projective factors
in the given CICY manifold. Then the harmonic $\left( 1,1\right) -$forms on $%
Y$ are simply the pullbacks of the K\"{a}hler forms $J_{1},...,J_{m}$ on the
corresponding complex projective factors. Now suppose we have some
involutions acting on $Y\times S^{1}$. First let us consider the case when
there are no $C$ involutions. In this case, no projective factors are mixed,
and each of the K\"{a}hler forms is odd under the involution. Hence in this
case, $h_{1,1}^{-}=h_{1,1}$ and $h_{1,1}^{+}=0$. From (\ref{barelybetti}),
we thus have on the $7$-dimensional quotient space that $b_{2}=0$ and $%
b_{3}=h_{11}+h_{2,1}+1.$

Now consider the case when we have one $C$ involution acting on $Y$. Without
loss of generality assume that the $C$ involution acts on the first two
projective factors. Then $J_{1}+J_{2}$ is odd, while $J_{1}-J_{2}$ is even
under this involution. The remaining K\"{a}hler forms remain odd as before.
So in this case, $h_{1,1}^{-}=h_{1,1}-1$ and $h_{1,1}^{+}=1$, and so $%
b_{2}=1 $ and $b_{3}=h_{1,11}+h_{2,1}$. When we have multiple $C$
involutions, $b_{2} $ correspondingly is equal to the number of $C$
involutions: 
\begin{subequations}
\label{cinvolbettis}
\begin{eqnarray}
b_{2} &=&n_{c} \\
b_{3} &=&h_{1,1}+h_{2,1}+1-n_{c}
\end{eqnarray}%
where $n_{C}$ is the number of $C$ involutions acting on the base CICY
manifold.

Thus far we have assumed that on the CICY manifold $h_{1,1}=m$. However this
is not always the case - in the list of CICY manifolds by Green et al, 4874
configurations satisfy this criterion, while the rest do not. The class of
CICYs for which this equality holds has been referred to as \emph{favourable 
}by Candelas and He \cite{CandelasHe}. It is known however, that there are
various identities which link together configuration matrices. One of the
simplest identities \cite{CandelasCICY1} is 
\end{subequations}
\begin{equation}
\left[ 
\begin{array}{c}
1 \\ 
X%
\end{array}%
\right\Vert \left. 
\begin{array}{c}
a+b \\ 
M%
\end{array}%
\right] =\left[ 
\begin{array}{c}
1 \\ 
1 \\ 
X%
\end{array}%
\right\Vert \left. 
\begin{array}{cc}
1 & a \\ 
1 & b \\ 
0 & M%
\end{array}%
\right]  \label{configid}
\end{equation}%
This is derived from the basic identity 
\begin{equation}
\left[ 
\begin{array}{c}
1 \\ 
1%
\end{array}%
\right\Vert \left. 
\begin{array}{c}
1 \\ 
1%
\end{array}%
\right] =\mathbb{CP}^{1}
\end{equation}%
which essentially says that a homogeneous hypersurface of degree $1$ in $%
\mathbb{CP}^{1}\times \mathbb{CP}^{1}$ is again $\mathbb{CP}^{1}$. Using (%
\ref{configid}) we can expand any configuration matrix which has a $\mathbb{%
CP}^{1}$ to an arbitrary size. In particular, if $h_{1,1}>m$ for the
original configuration matrix, we can expand the matrix so that it has
precisely $h_{1,1}$ projective factors. Once we have such a matrix, we again
find the possible involutions and calculate the Betti numbers of the
corresponding barely $G_{2}$ manifolds. Employing this procedure, we can
cover all but $37$ configurations.

After doing all the calculations we find the following pairs of Betti
numbers of the barely $G_{2}$ manifolds%
\begin{equation}
\begin{tabular}{ll}
$b_{2}$ & $b_{3}$ \\ 
$0$ & $31+2k\ $for $k=0,...,22,24,29,30$ \\ 
$1$ & $30+2k\ $for $k=0,...,19,21$ \\ 
$2$ & $29+2k\ $for $k=0,...,10,12,13,15$ \\ 
$3$ & $28+2k$ \ for $k=0,...,7,9,10$ \\ 
$4$ & \thinspace $27+2k$ $\ $for $k=0,...,3,5,7$ \\ 
$5$ & $26+2k$ \ for $k=0,1,3,4$ \\ 
$6$ & $25,31$ \\ 
$7$ & $24$%
\end{tabular}%
\end{equation}%
Thus we have a total of $84$ distinct pairs of Betti numbers. All of these
pairs have odd $b_{2}+b_{3}$ , and while most of Joyce's examples of $G_{2}$
holonomy manifolds have $b_{2}+b_{3}\equiv 3\ \func{mod}\ 4$, here we have a
mix between $b_{2}+b_{3}\equiv 1\ \func{mod}\ 4$ and $b_{2}+b_{3}\equiv 3\ 
\func{mod}\ 4$.

\section{Concluding remarks}

We have obtained the Betti numbers of barely $G_{2}$ manifolds obtained from
Complete Intersection Calabi-Yau manifolds. This gives a class of manifolds
that have an explicit description. One of the ways to use these examples is
to try and understand the moduli spaces. On one hand we know the structure
of the moduli space of the underlying CICY manifolds, but on the other hand,
previous general results about the structure of $G_{2}$ moduli spaces \cite%
{karigiannis-2007a,GrigorianYau1} could be applied to these specific cases.
In particular, quantities like the Yukawa couplings and curvature could be
calculated for these examples. This should then give a relationship between
the corresponding Calabi-Yau quantities and the $G_{2}$ quantities. This
could then lead to much better understanding of $G_{2}$ moduli spaces and
their relationship to Calabi-Yau moduli spaces.

Another direction could be to construct barely $G_{2}$ manifolds from some
larger class of Calabi-Yau manifolds. In particular it is interesting to see
what is the relationship between manifolds constructed from Calabi-Yau
mirror pairs, and whether this could shed some light on possible $G_{2}$
mirror symmetry.

\setcounter{equation}{0}

\bibliographystyle{jhep2}
\bibliography{refs2}

\end{document}